\documentclass{amsart}
\usepackage{amssymb,amsfonts,amsmath}
\usepackage[all]{xy}
\usepackage[shortlabels]{enumitem}
\usepackage{hyperref, color}		

\newtheorem{theorem}{Theorem}[section]

\newtheorem{corollary}{Corollary}[section]
\theoremstyle{definition}

\newtheorem{lemma}{Lemma}[section]
\newtheorem{question}{Question}[section]
\theoremstyle{remark}

\theoremstyle{case}

\theoremstyle{claim}

\theoremstyle{theoremA}

\theoremstyle{theoremB}
  
\newcommand{\real}{\mathbb{R}}
\newcommand{\cpx}{\mathbb{C}}

\newcommand{\hy}{\mathbb{H}}
\newcommand{\Sp}{\mathbb{S}}

\newcommand{\al}{\alpha}

\newcommand{\ep}{\epsilon}

\newcommand{\be}{\beta}
\newcommand{\ka}{\kappa}

\newcommand{\de}{\delta}

\newcommand{\weil}{\mathcal{W}}
\newcommand{\dd}{\,\textrm{d}V}

\newcommand{\lan}{\left\langle}
\newcommand{\ran}{\right\rangle}
\newcommand{\p}{\partial}

\newcommand{\m}{\mathcal}
\newcommand{\tr}{\mathrm{tr}\,}

\newcommand{\mk}{\mathfrak}

\newcommand{\vp}{\varphi}

\newcommand{\setap}{\rightarrow}
\begin{document}

\title[Einstein $4-$Manifolds and Nonpositive Curvature]{Einstein $4-$Manifolds and Nonpositive Isotropic Curvature}
%

\author{A. Brasil Jr}

\author{E. Costa}

\author{F. Vit\'orio}
\date{}
\subjclass[2010]{Primary 53C25; Secondary 53C24}
\begin{abstract} This note is devoted to study the implications of nonpositive isotropic curvature and  negative Ricci curvature for Einstein $4-$manifolds.  \end{abstract}
\maketitle
\thispagestyle{empty}
\section{Introduction}

Let $M$ be an oriented $4-$dimensional Riemannian manifold. For each point of $M$, we can consider the complexification $T_pM \otimes \mathbb{C}$ of the tangent space $T_pM$ at the point $p$. There are two natural extensions of the inner product on $T_pM$ to the complexified tangent space  $T_pM \otimes \mathbb{C}$, namely, as a complex bilinear form $(,)$ or as a Hermitian inner product $\langle, \rangle$. These extensions are related by $\langle u,v \rangle=(u,\bar v)$ for $u,v \in T_pM \otimes \mathbb{C}$. The tensor curvature of $M$ induces the  curvature operator \[\mathcal{R}: \Lambda^2T_pM\to  \Lambda^2T_pM\] on the fiber bundle of $2-$forms $\Lambda^2 M$ . After complexification, it extends to complex vectors by linearity, and we can define the complex sectional curvature of a two-dimensional subspace $\pi$ of $T_pM \otimes \mathbb{C}$ by
\[
K(\pi)= \langle \mathcal{R}(u\wedge v), \bar v \wedge u \rangle,
\]
where $\{u,v\}$ is any unitary basis of $\pi$. A vector $u\in T_pM \otimes \mathbb{C}$ is said isotropic if $(u,u)=0.$ While a subspace $E \subset T_pM \otimes \mathbb{C}$ is said isotropic if every vector $u\in E$ is isotropic. Moreover, $M$ is said to have \emph{nonpositive isotropic curvature} if $K(\pi)\leq 0$ for all isotropic two-dimensional subspace $\pi$ of $T_pM \otimes \mathbb{C}$. The notions of \emph{nonnegative}, \emph{negative} and \emph{positive isotropic curvature} can be defined in a similar way.

 It is well known that for an oriented $4$-dimensional Riemannian manifold $M^4$ the Hodge star operator splits the fiber bundle of $2-$forms $\Lambda^2 M= \Lambda^+ \oplus \Lambda^-$, where $\Lambda^\pm=\left\{ \omega \in \Lambda^2 M\big|*\omega=\pm \,\omega \right\}$ denote the $\pm$eigenspaces of that operator. At this context, the Weyl curvature tensor $\weil$ comutes with the Hodge star operator, and it is, therefore, an endomorphism of $\Lambda^2 M$, such that $\weil= \weil^+ \oplus \weil^-$.  Using such a decomposition,  we see that the non positivity of the isotropic  curvature operator is equivalent to the non positivity of the endomorphisms $P^\pm= \frac{s}{6}I_{\Lambda^+} \pm \weil^\pm,$ where $s$ stands for the scalar curvature of $M^4.$

In \cite{MW}, Micallef and Wang proved that
if $M$ is a compact Einstein $4$-manifold with nonnegative isotropic curvature,
then $M$ is locally symmetric. Recently, Brendle \cite{brendle} proved that a compact Einstein manifold with nonnegative isotropic curvature must be a locally symmetric space of compact type. On the other hand, Seshadri \cite{Se} showed that
any compact Riemannian 4-manifold $M$ admits a Riemmanian metric with
strictly negative isotropic curvature. Despite the result of Seshadri, we can
propose the following questions:

\begin{question} \label{question1} Is a compact Einstein 4-manifold with nonpositive
isotropic curvature locally symmetric?
\end{question}

Since nonpositive curvature operator implies nonpositive isotropic curvature,
we can  propose the following weak version of Question \ref{question1}.

\begin{question} \label{question2} Is a compact Einstein 4-manifold with nonpositive
curvature operator locally symmetric?
\end{question}

Very recently, Fine, Kasnov and Panov \cite{FKP} have proposed a similar question. More precisely, they posed the following question.

\begin{question} Do there exists compact Einstein $4-$manifolds with
scalar curvature $s<0$, for which $\frac{s}{12}I+\weil^+$ is negative definite,
besides hyperbolic and complex-hyperbolic $4-$manifolds?
\end{question}

It is important to highlight that if an oriented $4-$manifold has endomorphism $\frac{s}{12}+\weil^+$ negative definite, then the endomorphism $P^+=\frac{s}{6}-\weil^+$ is also negative definite.

There are some evidences in favor of a positive answers to Questions \ref{question1}
and \ref{question2}. More precisely, we have the following useful informations:

\begin{enumerate}

\item[(a)] The only known examples of compact Einstein $4-$manifolds $M$ with negative Ricci curvature and nonpositive isotropic curvature are locally
symmetric spaces whose universal covering of $M$ is isometric to either complex hyperbolic space $\mathbb{CH}^2,$ real hyperbolic space $\mathbb{H}^4_c$ or a product
of two hyperbolic spaces $\mathbb{H}^2_c \times \mathbb{H}^2_c.$\\


\item[(b)] $4-$manifolds of nonpositive curvature operator or negative $1/4-$pinched sectional curvature have nonpositive isotropic curvature.\\

\item[(c)] In \cite{Vi}, Ville proved that if a compact oriented Riemannian $4-$manifold
M has negative $1/4-$pinched sectional curvature, then the Euler characteristic $\chi$ and its signature $\tau$ satisfy $\chi\geq 3|\tau|.$ Moreover, the equality occurs if and only if the universal covering of $M$ is isometric to $\mathbb{CH}^2$.\\

\item[(d)] In \cite{ZY}, Zheng and Yau proved  that if a compact K\"ahler-surface has negative $1/4-$pinched sectional curvature, then the universal covering
of $M$ is isometric to complex hyperbolic space $\mathbb{CH}^2$ with its standard metric.\\

\item[(e)] In \cite{Zh}, F. Zheng showed that if a compact K\"ahler-surface $M$ has nonpositive sectional curvature, then $M$ has signature $\tau \geq 0.$ Moreover, $\tau=0$ if and only if the universal covering of $M$ is isometric to $\mathbb{H}^2_c \times \mathbb{H}^2_c.$
\end{enumerate}

Proceeding, it is well-known that a compact oriented Einstein $4-$manifold satisfies
\begin{equation}\label{eulerformula}
8\pi^2 \chi = \int_M \big(|\weil|^2+\frac{s^2}{24} \big) \dd,
\end{equation}
and
\begin{equation}\label{signatureformula}
12 \pi^2 \tau = \int_M \big(|\weil^+|^2-|\weil^-|^2 \big) \dd.
\end{equation} These formulae tell us that a compact oriented Einstein $4$-manifold must to satisfy the well-known Hitchin-Thorpe inequality:
\[
\chi \geq \frac{3|\tau|}{2}.
\]

Our first result gives a similar obstruction to the existence of Einstein metrics with nonpositive isotropic curvature on $4$-manifolds. More precisely, we have the following result.

\begin{theorem} \label{thm1} Let $M$ be a compact oriented Einstein $4-$manifold with negative Ricci curvature  and nonpositive isotropic curvature. Then we have:
\begin{equation} \label{hitchinthorpe} \chi \geq \frac{15}{8}|\tau|.
\end{equation} In addition, if $M$ is a K\"ahler-Einstein manifold, then
\begin{equation} \label{einstein-Kaehler1}
 \chi \geq 3 |\tau|.
 \end{equation} Moreover, the equality in (\ref{einstein-Kaehler1}) occurs if and only if the universal covering $\widetilde{M}$ of $M$ is isometric to $\mathbb{CH}^2$.
\end{theorem}

\begin{theorem} \label{thm2} Let $M$ be a compact oriented Einstein $4-$manifold with negative Ricci curvature and volume $V$.
\begin{enumerate}
\item If $M$ has nonpositive curvature operator, then $\chi \geq 3 |\tau|$. Moreover, the equality occurs if and only if the universal covering $\widetilde{M}$ of $M$ is isometric to $\mathbb{CH}^2.$
\item If $M$ has nonpositive curvature operator, then $\chi \leq \frac{\rho^2 V}{4\pi^2}$. Moreover, the equality occurs if and only if the universal covering  $\widetilde{M}$ of $M$ is isometric to $\mathbb{H}^2_c \times \mathbb{H}^2_c.$
\item If the sectional curvature $K$ of $M$ satisfies $\sup K \leq \frac{\rho}{6},$ then $\chi \leq \frac{\rho^2 V}{6\pi^2}$. Furthermore, the equality occurs if and only if the universal covering $\widetilde{M}$ of $M$ is isometric to $\mathbb{CH}^2.$
\end{enumerate}
\end{theorem}


\section{Proof of Theorem \ref{thm1}}
The first statement  of Theorem \ref{thm1} is a straightforward consequence of the following lemma.

\begin{lemma}\label{lemma1}Let $M$ be a compact oriented Einstein $4-$manifold with negative Ricci curvature. If $M$ has nonpositive isotropic curvature, then
\[
2\chi+3|\tau| \leq \frac{3\rho^2 V}{2\pi^2}\leq 18 \chi-27|\tau|.
\]
\end{lemma}
\noindent \textbf{Proof:}
We follow the ideas developed in \cite{Co}. Indeed, for each $x\in M$ there exist an orthonormal basis of $T_xM$ and a corresponding orthonormal basis of $\Lambda^\pm$, such that $\weil^\pm$ has the respective eigenvalues
\[
\lambda_1\pm \mu_1-\rho/3 \leq \lambda_2\pm \mu_2-\rho/3 \leq \lambda_3\pm \mu_3-\rho/3,
\]
where $\lambda_1, \lambda_2,$ and $\lambda_3$ are the principal sectional curvatures of $M.$ Hence, $\lambda_1+\lambda_2+\lambda_3=\rho$ and $\mu_1+\mu_2+\mu_3=0$.  It will be convenient to see $\lambda=(\lambda_1, \lambda_2, \lambda_3)$ and $\mu=(\mu_1, \mu_2, \mu_3)$ as vectors in $\mathbb{R}^3$ with canonical inner product $\langle , \rangle$ and its canonical norm $|.|$.

Using these settings,  (\ref{eulerformula}) and (\ref{signatureformula}) become
\begin{equation}\label{eulerformula2}
4\pi^2 \chi = \int_M   \big(|\lambda|^2+|\mu|^2  \big)\dd,
\end{equation}
and
\begin{equation}\label{signatureformula2}
3 \pi^2 \tau = \int_M \langle \lambda , \mu  \rangle\dd.
\end{equation}

We point out that the non positivity of the isotropic curvature of $M$:
\[ R_{1313}+R_{1414}+R_{2323}+R_{2424}+2R_{1234} \leq 0
\] is equivalent to $$\lambda_1\pm\mu_1-\rho \geq 0.$$
Then, by setting $\alpha_i^{\pm}=\lambda_i \pm \mu_i -\rho,$ we have
\begin{equation} \label{monotonicity}
0\leq\alpha_1^{\pm} \leq \alpha_2^{\pm} \leq \alpha_3^{\pm}.
\end{equation}
Now, we notice that (\ref{monotonicity}) implies
\[ 4\rho^2= (\alpha_1^{\pm} + \alpha_2^{\pm} + \alpha_3^{\pm})^2\geq (\alpha_1^{\pm})^2 + (\alpha_2^{\pm})^2 + (\alpha_3^{\pm})^2.
\]
But
\[
(\alpha_1^{\pm})^2 + (\alpha_2^{\pm})^2 + (\alpha_3^{\pm})^2= |\lambda|^2+|\mu|^2\pm 2 \langle \lambda, \mu \rangle + \rho^2.
\]
From here it follows that
\[
 |\lambda|^2+|\mu|^2\pm 2 \langle \lambda, \mu \rangle \leq 3\rho^2.
\]
From (\ref{eulerformula2}) and (\ref{signatureformula2}) we deduce
\begin{equation} \label{forward}
2\chi+3|\tau| \leq \frac{3\rho^2V}{2 \pi^2}.
\end{equation}
On the other hand, using once more (\ref{eulerformula2}) and (\ref{signatureformula2}) we arrive at
\[
\begin{array}{rcl}
2\pi^2\big( 2\chi\pm 3\tau\big)&= &\displaystyle{\int_M \left( |\lambda|^2+|\mu|^2 \pm 2 \langle \lambda, \mu \rangle\right) \dd}\\ \\
&=&\displaystyle{\int_M |\lambda\pm \mu|^2 \dd}\geq \frac{\rho^2 V}{3},
\end{array}
\]
so that
\begin{equation} \label{reverse}
2\pi^2\big( 2 \chi -3 |\tau|\big) \geq  \frac{\rho^2 V}{3}.
\end{equation}
Finally, by combining (\ref{forward}) and (\ref{reverse}) we get
\[
2\chi+3|\tau| \leq \frac{3\rho^2 V}{2\pi^2}\leq 18 \chi-27|\tau|.
\]
This concludes the proof of lemma.
\smallskip

In order to prove the last statement of Theorem \ref{thm1} we consider $M$ to be a K\"ahler-Einstein manifold with negative Ricci curvature $\rho$. In this case, the eigenvalues of $\weil^+$ are $2\rho/3$, $-\rho/3$ and $-\rho/3$. Thus, we can use (\ref{eulerformula}) and (\ref{signatureformula}) to infer
\begin{equation} \label{remain1}
8\pi^2 \chi = \int_M \big( |\weil^-|^2+4\rho^2/3 \big) \dd
\end{equation}
and
\begin{equation} \label{remain2}
12\pi^2\tau= \int_M \big( 2\rho^2/3 - |\weil^-|^2\big)\dd.
\end{equation}
These two above inequalities gives
\[
8\pi^2 \big(\chi -3 \tau \big)=3  \int_M  |\weil^-|^2\dd.
\]
So, we conclude that $\chi \geq 3 \tau$.

On the other hand, since that $M$ has nonpositive isotropic curvature, we deduce that the eigenvalues of $\weil^-$ satisfy $\gamma_3\geq \gamma_2 \geq \gamma_1 \geq 2\rho/3$ and $\gamma_1+\gamma_2+\gamma_3=0$. From here it follows that $0\leq \gamma_3\leq -4\rho/3$ and $2\rho/3 \leq \gamma_1\leq 0$. With this setting, a simple computation yields
\begin{equation}\label{remainweil}
|\weil^-|^2= \gamma_1^2+\gamma_2^2+\gamma_3^2 \leq 8\rho^2/3.
\end{equation}
Sum (\ref{remain1}) and (\ref{remain2}), and then we plug (\ref{remainweil}) to obtain $\chi\geq -3\tau.$ Hence, $\chi\geq 3|\tau|$.

Moreover, if $\chi=3|\tau|$, then $\weil^-=0$. Therefore, the universal covering of $M$ is isometric to $\mathbb{CH}^2$. This finishes the proof of Theorem \ref{thm1}.

\section{Proof of Theorem \ref{thm2}}

The first and second assertions follow directly from the following

\begin{lemma}\label{lemma2}Let $M$ be a compact oriented Einstein $4-$manifold with negative Ricci curvature.
 If $M$ has nonpositive  curvature operator, then
\[
2\chi+3|\tau| \leq \frac{\rho^2 V}{2\pi^2}\leq 6 \chi-9|\tau|.
\]
\end{lemma}

\noindent \textbf{Proof:} First of all, observe that the eigenvalues of the curvature operator of $M$ are $\lambda_i\pm \mu_i$, for $i=1,2,3$. Hence, the non positivity of the curvature operator implies that $\lambda_i\pm \mu_i \leq 0$. Using this data, we have
\[
\begin{array}{rcl}
\rho^2=\left(\sum_i(\lambda_i\pm\mu_i)\right)^2&=& |\lambda\pm\mu|^2+2\sum_{i<j}(\lambda_i\pm\mu_i)
(\lambda_j\pm\mu_j)\\ \\
&\geq& |\lambda \pm \mu|^2=|\lambda|^2+|\mu|^2 \pm 2 \langle \lambda, \mu \rangle,
\end{array}
\]
Thus
\begin{equation} \label{estimatel2}
|\lambda|^2+|\mu|^2  \pm 2 \langle \lambda, \mu \rangle \leq \rho^2.
\end{equation}
Upon integrating (\ref{estimatel2}) we use  (\ref{eulerformula2}) and (\ref{signatureformula2}) to arrive at
\begin{equation}\label{estimatel22}
2\chi+3|\tau| \leq \frac{\rho^2V}{2\pi^2}.
\end{equation}
Clearly, the combination of (\ref{estimatel22}) and (\ref{reverse}) yields
$$
2\chi+3|\tau| \leq \frac{\rho^2 V}{2\pi^2}\leq 6 \chi-9|\tau|.
$$ This concludes the proof of  Lemma \ref{lemma2}.

From now on we shall characterize the equality in the first and second assertions of Theorem \ref{thm2}. Initially, we notice that if $\chi=3|\tau|$, then the universal covering of $M$ is  isometric to $\mathbb{CH}^2.$ Next, if $\chi=\frac{\rho^2 V}{4\pi^2}$, then we have from (\ref{estimatel22}) that $\tau=0$. This implies that $\lambda_i\pm \mu_i=0$, for $i=2,3$. So, the eigenvalues of
the curvature operator of $M$ are constant and $M$ is locally symmetric.  Therefore, $\weil^+ \not \equiv 0$ and $\tau=0,$ which tell us that $\widetilde{M}$ is isometric to $\mathbb{H}^2_c \times \mathbb{H}^2_c$.

In order to proof the third item of Theorem \ref{thm2} we suppose that the sectional curvature $K$ of $ M$ satisfies $\sup K \leq \rho/6.$ We then use the inequality (2.8) in \cite{Co} to infer  $\chi \leq \frac{\rho^2 V}{6\pi^2}$. Moreover, if  $\chi = \frac{\rho^2 V}{6\pi^2}$, then is locally symmetric. Therefore, $\widetilde{M}$ is isometric to $\mathbb{CH}^2,$ as we wanted to prove.

\vspace{0.7cm}
\begin{small}
\begin{tabular}{lcr}
\begin{tabular}{l}
Aldir Brasil J\'unior\\
Universidade Federal do Cear\'a\\
Departamento de Matem\'{a}tica\\   60.455-760, Fortaleza, CE
\\
Brazil\\
\verb+aldir@mat.ufc.br+\\
\end{tabular}\

&\quad \quad&

\begin{tabular}{l}
\'Ezio Costa\\
Universidade Federal da Bahia\\
Instituto de Matem\'{a}tica\\
40170-110 Salvador-BA\\
Brazil\\
\verb+ezio@ufba.br+\\
\end{tabular}\\ \\
\begin{tabular}{l}
Feliciano Vit\'orio\\
Universidade Federal de Alagoas\\
Instituto de Matem\'{a}tica\\
57072-900 Macei\'o-AL\\
Brazil\\
\verb+feliciano@pos.mat.ufal.br+\\
\end{tabular}\

\end{tabular}

\vspace{1cm}

\end{small}

\end{document}